\ifdefined\singlecolumn
\documentclass[11pt,draftclsnofoot,onecolumn]{IEEEtran}
\else
\documentclass[10pt,twocolumn]{IEEEtran}
\fi

\usepackage{cite}
\usepackage{graphicx} 
\graphicspath{{figures/}}

\usepackage[tbtags]{amsmath}
\usepackage{amssymb,color,amsthm}
\usepackage[none]{hyphenat}
\usepackage{subfigure}
\usepackage{bbm}
\usepackage{etoolbox}
\usepackage[dvipsnames]{xcolor}
\usepackage{url}
\usepackage{flushend}
\newtoggle{twocolumn}
\ifCLASSOPTIONtwocolumn\toggletrue{twocolumn}\else\togglefalse{twocolumn}\fi

\def\IS{\text{IS}}
\def\KL{\text{KL}}

\def\minimize{\mathop{\text{minimize}}}

\usepackage[utf8]{inputenc}

\usepackage{float}
\restylefloat{table}

\usepackage{caption}

\usepackage[ruled,noline]{algorithm2e}
\usepackage{multirow}
\usepackage{textcomp}
\let\labelindent\relax
\usepackage{enumitem}

\newtheorem{lemma}{Lemma}
\newtheorem{theorem}{Theorem}
\newtheorem{corollary}{Corollary}
\newtheorem{define}{Definition}
\newtheorem{remark}{Remark}
\newtheorem{prop}{Proposition}
\newtheorem{assump}{Assumption}

\newcommand{\bigo}[1]{\mathcal{O}(#1)}
\newcommand{\eg}{{\it e.g.}}
\newcommand{\ie}{{\it i.e.}}
\newcommand{\norm}[1]{\|#1\|}
\newcommand{\argmin}{\text{argmin}}

\newcommand{\reals}{{\mbox{$\mathbf{R}$}}}

\newcommand{\size}[1]{\left|#1\right|}
\newcommand{\abs}[1]{\size{#1}}
\newcommand{\inn}[2]{\langle#1, #2\rangle}

\def\a{{\mathbf a}}
\def\b{{\mathbf b}}
\def\x{{\mathbf x}}
\def\y{{\mathbf y}}
\def\z{{\mathbf z}}
\def\u{{\mathbf u}}
\def\v{{\mathbf v}}
\def\w{{\mathbf w}}
\def\d{{\mathbf d}}

\def\E{{\mathbb E}}

\def\A{{\mathbf A}}
\def\X{{\mathbf X}}

\def\cX{{\mathcal X}}

\begin{document}

\title{Randomized Bregman Coordinate Descent Methods for Non-Lipschitz Optimization}
\author{Tianxiang~Gao, \emph{Student Member, IEEE},
	Songtao~Lu, \emph{Member, IEEE},
	Jia~Liu, \emph{Senior Member, IEEE},
    \\and Chris~Chu, \emph{Fellow, IEEE}


\thanks{
    Tianxiang Gao and Chris Chu are with the Department of Electrical and Computer Engineering,
    Iowa State University, Ames, IA 50011, USA (emails: \{gaotx,cnch\}@iastate.edu).}
\thanks{
    Songtao Lu is with IBM Research AI, IBM Thomas J. Waston Research Center, Yorktown Heights, New York 10562,  USA (email: songtao@ibm.com). }
\thanks{
	Jia Liu is with the Department of Computer Science,
	Iowa State University, Ames, IA 50011, USA (email: jialiu@iastate.edu).}
}

\maketitle

\begin{abstract}
    We propose a new \textit{randomized Bregman (block) coordinate descent} (RBCD) method for minimizing a composite problem, where the objective function could be either convex or nonconvex, and the smooth part are freed from the global Lipschitz-continuous (partial) gradient assumption. Under the notion of relative smoothness based on the Bregman distance, we prove that every limit point of the generated sequence is a stationary point. Further, we show that the iteration complexity of the proposed method is $\bigo{n\varepsilon^{-2}}$ to achieve $\epsilon$-stationary point, where $n$ is the number of blocks of coordinates. If the objective is assumed to be convex, the iteration complexity is improved to $\bigo{n\epsilon^{-1}}$. If, in addition, the objective is strongly convex (relative to the reference function), the global linear convergence rate is recovered. We also present the accelerated version of the RBCD method, which attains an $\bigo{n\varepsilon^{-1/\gamma} } $ iteration complexity for the convex case, where the scalar $\gamma\in [1,2]$ is determined by the \textit{generalized translation variant} of the Bregman distance. Convergence analysis without assuming the global Lipschitz-continuous (partial) gradient sets our results apart from the existing works in the composite problems.

\end{abstract}

\section{Introduction}
\label{sec:intro}
In this paper, we consider a composite optimization problem in the following form
\begin{align}
\minimize_{\x}\; F(\x)\equiv f(\x) + r(\x),\label{opt:problem}
\end{align}
where $r$ has $n$ separated blocks. More specifically, we have
\begin{align}
r(\x) = \sum_{i=1}^{n} r_i(\x_i),
\end{align}
where $\x_i$ denotes a subvector of $\x$ with dimension $N_i$ such that $\sum_{i=1}^n N_i = N$, and each $r_i$ is a (possibly nonsmooth) convex function.

Due to the block separable structure, Problem \eqref{opt:problem} can be solved by \textit{(block) coordinate descent} (CD) methods and/or their variants, especially in the large scale optimization problems. Roughly speaking, these methods are based on the strategy of selecting one coordinate/block of variables at each iteration using some index selection procedure (\eg, cyclic, greedy, randomized). This often dramatically reduces the computational complexity of the algorithms per iteration as well as memory storage, making these methods simple and salable. See for instance \cite{nesterov2012efficiency,lu2015complexity,nutini2015coordinate,lu2019pa,beck2013convergence} and references therein and a short summary in Table~\ref{tab:bcd summary}, as well as the recent comprehensive review paper \cite{wright2015coordinate} for the up-to-date materials.

\begin{table*}
	\centering
	\begin{tabular}{ |c | c | c |c| c| c|}
		\hline
		Paper & Algorithm & Problem& Lipschitz-continuous (partial) gradient & Iteration complexity & Acceleration\\  \hline
		\cite{nesterov2012efficiency,lu2015complexity}& Randomized CD & Convex & Needed& \multirow{3}{*}{$\bigo{n\varepsilon^{-1}}$} &No\\ \cline{1-4} \cline{6-6}
		\cite{beck2013convergence,saha2013nonasymptotic,hong2015unified,sun2015improved}& Cyclic CD & Convex & Needed&  &No\\ \cline{1-4} \cline{6-6}
		\cite{nesterov2012efficiency,dhillon2011nearest,nutini2015coordinate}& Greedy CD & Convex & Needed&  &No\\ \hline
		\hline
		\cite{nesterov2012efficiency,lu2015complexity,lin2014accelerated,fercoq2015accelerated,qu2016coordinate}& Randomized CD & Convex & Needed& \multirow{2}{*}{$\bigo{n\varepsilon^{-1/2}}$} &Yes\\ \cline{1-4} \cline{6-6}
		\cite{song2017accelerated,lu2018accelerating,locatello2018matching}& Greedy CD & Convex & Needed& &Yes\\ \hline
		\hline
		\cite{patrascu2015efficient}& Randomized CD & Nonconvex & Needed& \multirow{2}{*}{$\bigo{n\varepsilon^{-2}}$} &No\\  \cline{1-4} \cline{6-6}
		\cite{shi2017inexact}& Cyclic CD & Nonconvex & Needed& &No\\ \hline
		\hline
		\cite{birnbaum2011distributed,bauschke2016descent,lu2018relatively}& GD& Convex & No needed& $\bigo{\varepsilon^{-1}}$ &No\\ \hline
		\cite{hanzely2018accelerated}& GD& Convex & No needed& $\bigo{\varepsilon^{-1/\gamma}}$ &Yes\\ \hline
		\cite{bolte2018first}& GD& Nonconvex & No needed& $\bigo{\varepsilon^{-2}}$ &No\\ \hline
		\cite{hanzely2018fastest}& Randomized CD& Convex & No needed& $\bigo{n\varepsilon^{-1}}$ &No\\ \hline
		\cite{gao2019leveraging}& Greedy CD& Nonconvex & No needed& $\bigo{n\varepsilon^{-2}}$ &No\\ \hline
		\textbf{This paper}& Randomized CD& Convex & No needed& $\bigo{n\varepsilon^{-1}}$ &No\\ \hline
		\textbf{This paper}& Randomized CD& Convex & No needed& $\bigo{n\varepsilon^{-1/\gamma}}$ &Yes\\ \hline
		\textbf{This paper}& Randomized CD& Nonconvex & No needed& $\bigo{n\varepsilon^{-2}}$ &No\\ \hline
		\hline
	\end{tabular}
	\caption{Summary of reference function and relative smoothness coefficient for various $\beta$ values.}
	\label{tab:bcd summary}
\end{table*}

A widely used assumption in showing the convergence of CD methods in the literature is that the (partial) gradient of $f$ is globally Lipschitz-continuous. However, this could be a restrictive assumption violated in diverse applications in practice, such as matrix factorization \cite{lee1999learning}, tensor decomposition \cite{kim2007nonnegative}, matrix/tensor completion \cite{xu2012alternating}, Poisson likelihood models \cite{he2016fast}, etc. Although this assumption may be relaxed by adopting conventional line search methods, the efficiency and computational complexity of the first-order method are unavoidably distorted, especially when the size of the problem is large. In fact, this longstanding issue also appears in the classical \textit{proximal gradient descent} (PGD) method. Fortunately, this issue is solved in \cite{birnbaum2011distributed,bauschke2016descent,lu2018relatively}. They develop a new framework called \textit{Bregman proximal gradient} (BPG) method that adapts the geometry of $f$ by the Bregman distance. In such a way,  the decrease of the objective value can be still quantified. As a result, they are able to characterize the convergence behavior of BPG for minimizing convex composite problems without assuming globally Lipschitz-continuous gradient of the objective function. Further, this framework has been extended  to the case of nonconvex optimization in \cite{bolte2018first}.

Despite the crucial issue is solved in PGD-type methods, there are only few results on CD-type methods. A \textit{cyclic Bregman coordinate descent} (CBCD) method has been proposed in \cite{ahookhosh2019multi,wang2018block}, but no rates are given. In \cite{hanzely2018fastest}, the authors provide the convergence rate result using randomized (block) coordinate selection strategy in a special case where $F$ is smooth convex and $r\equiv 0$. To the best of our knowledge, how to deal with this crucial issue is still an open problem, when using CD methods to solve a nonsmooth and convex/nonconvex Problem \eqref{opt:problem}. Furthermore, the accelerated version of the RBCD method has not been proposed yet, and its iteration complexity analysis is still  open as well. In this paper, we bridge these gaps by proposing a \textit{randomized Bregman (block) coordinate descent} (RBCD) method and its accelerated variant. The comprehensive convergence analyses are established. The main contributions are highlighted as follows.
\begin{enumerate}[leftmargin=4mm]
	\item We propose a randomized Bregman (block) coordinate descent (RBCD) method to solve the composite problem where the smooth part does not have the global Lipschitz-continuous (partial) gradient property.
	\item By adapting the relative smoothness framework, we establish a rigorous convergence rate analysis of the RBCD method, showing that the convergence rate to an stationary point is $\bigo{n\varepsilon^{-2}}$ if $F$ is nonconvex, where $k$ is the number of iterations.
	\item If $F$ is convex, RBCD achieves the global sublinear convergence rate of $\bigo{n\varepsilon^{-1}}$. The global linear convergence rate is obtained if $f$ is (relative) strongly convex.
	\item The RBCD method can also be accelerated in the relative smoothness setting. The iteration complexity of $\bigo{n\varepsilon^{-1/\gamma} }$ can be obtained through the notion of \textit{generalized translation variant} (explained in the latter section) of the Bregman distance.
\end{enumerate}

\section{Preliminaries}
\textbf{Notation}. Throughout this paper, we use bold upper case letters denote matrices (\eg. $\X$), bold lower case letters denote vectors (\eg, $\x$), and Calligraphic letters (\eg, $\mathcal{X}$) are used to denote sets. We use $\norm{\cdot}$ to denote the Euclidean norm. $\delta_{\mathcal{X}}(\x)$ represents the indicator function: $ \delta_{\mathcal{X}}(\x) = 0$ if $\x\in \mathcal{X}$; otherwise, $ \delta_{\mathcal{X}}(\x) =\infty$. If $\mathcal{X}=\reals^N_+$, the indicator function becomes $\delta_+(\x)$. For a function $f$, $\nabla f(\x)$ denotes its the gradient, while $\nabla_i f(\x)$ is the partial gradient with respect to the $i$-th block. Let $f_i(\x_i)$ be the function with respect to the $i$-th block, while the rest of blocks are fixed. Clearly, we have $\nabla_i f(\x) = \nabla f_i (\x_i)$. If $f$ is not differentiable, $\partial f$ denotes the subdifferential of $f$.

Given a convex function $\phi$, the \textit{Bregman proximal mapping} of $\phi$ at a point $\x$ is defined as
\begin{align}
T_{\phi}(\x) = \argmin_\u \phi( \u) + D_h(\u, \x),
\end{align}
where $D_h(\u,\x) = h(\u) - h(\x) - \inn{\nabla h(\x)}{\u-\x}$ is the \textit{Bregman distance} with the reference convex function $h$. This mapping is well-defined since the functions $\phi$ and $h$ are convex. The convexity of $h$ also implies $D_h(\x,\y)\geq 0, \forall \x,\y$. If, in addition, $h$ is strictly convex, $D_h(\x,\y) = 0$ if and only if $\x=\y$. In the rest of this paper, we assume $h$ is strictly convex. Note that $D_h(\x,\y)$ is not symmetric in general. Therefore, we use \textit{symmetric coefficient} $\theta(h)$, defined by
\begin{align}
	\theta(h)=\inf_{\x \neq \y} \left\{ D_h(\x,\y)/D_h(\y,\x)  \right\},
\end{align}
to measure the symmetry. When $\phi=\delta_{\mathcal{\x}}$, the Bregman proximal mapping reduces to the \textit{Bregman projection}
\begin{align}
P_{\mathcal{X}}^h(\x)=\argmin\{D_h(\u,\x):\u\in \mathcal{X}\}.
\end{align}

\textbf{Problem Formulation}. Our goal is to solve the following composite optimization problem
\begin{align}
\minimize\; F(\x)\equiv f(\x) + r(\x),\label{problem}
\end{align}
where the following assumptions are made throughout this paper.
\begin{assump}\label{assume:basic}
	\item[(i)] $f$ is continuously differentiable.
	\item[(ii)] $r$ is convex, block separable, proper and loser semi-continuous.
	\item [(iii)] $F^* = \inf_\x F(\x) > -\infty$.
\end{assump}
An estimate $\x$ is said to be a \textit{stationary point} of $F$ if it satisfies
\begin{align}
0\in \partial F \equiv \nabla f(\x) + \partial r(\x).
\end{align}
Note that the objective function $F$ could be convex or nonconvex since we don't make the convexity assumption of $f$, which is the case in \cite{hanzely2018fastest}. In addition, the function $r$ could be an indicator function of a closed convex set, so that the problem formulation in \eqref{problem} includes the case where minimizing a nonsmooth objective function over a closed convex set.

\section{Randomized Bregman Coordinate Descent}
In this section, we introduce the \textit{randomized Bregman (block) coordinate descent} (RBCD) method for solving problem \eqref{problem}. Given the current estimate $\x$, the $i$-th block of coordinates is selected uniformly at random, then the new estimate $\x^+$ is updated as follows
\begin{align}
\x_i^+ = T_i(\x), \quad\text{and}\quad \x_j^+ = \x_j,\forall j\neq i,\label{eq:x^+}
\end{align}
where, for some stepsize $\alpha$, the vector $T_i(\x)$ is defined as
\begin{align}
	T_i(\x)=\argmin_{\u_i} \inn{\nabla_i f(\x)}{\u_i  - \x_i} + \frac{1}{\alpha}D_h(\u_i, \x_i) + r_i(\u_i).\label{eq:T_i}
\end{align}
Note that we drop the index $i$ in $D_{h_i}$ to simplify the notation. The algorithm is summarized in Algorithm~\ref{alg}.
\begin{algorithm}[h]
	\caption{Randomized Bregman (Block) Coordinate Descent (RBCD).}
	\label{alg}
	Choose $\x^0$.\\
	\For{$k=1,2\cdots$}{
		Choose $i_k\in \{1,2,\cdots, n\}$ uniformaly at random\\
		Specify the stepsize $\alpha^k$\\
		Compute $T_{i_k}(\x^k)$ from \eqref{eq:T_i}\\
		Update $\x^{k+1}$ by \eqref{eq:x^+}
	}
\end{algorithm}

Here the stepsize $\alpha$ can be determined by a conventional line search method and the global convergence results can be established. However, line search methods are usually expensive since this subroutine requires to evaluate the objective function multiple times to ensure the sufficient descent in the objective value. To establish convergence results for a CD-type method with a constant stepsize, the common assumption is that $\nabla f(\x)$ (or $\nabla_i f(\x)$) is globally Lipschitz-continuous \cite{nesterov2012efficiency,lu2015complexity,bonettini2016cyclic}. However, this assumption may be restrict to some modern optimization problems. See for instances \cite{lee1999learning,kim2007nonnegative,xu2012alternating,he2016fast} and reference therein. In the following section, we review the notion of \textit{relative smoothness} introduced in \cite{birnbaum2011distributed,bauschke2016descent,lu2018relatively}. This notion allows us to establish the convergence results for RBCD method without the assumption of global Lipschitz-continuous gradient.

\section{Convergence Analyses of RBCD}
We start with the definition of \textit{relative smoothness} \cite{lu2018relatively,bauschke2016descent}, by which a new descent lemma is obtained without the assumption of the global Lipschitz-continuity of (partial) gradient.
\begin{define}[Relative Smoothness]\cite[Definition~1.1]{lu2018relatively}
	A pair of functions $(g, h)$ are said to be relatively smooth if $h$ is convex and there exists a scalar $L>0$ such that $Lh-g$ is convex.
\end{define}
Moreover, the relative smoothness nicely translates the Bregman distance to produce a non-Lipschitz descent lemma \cite{lu2018relatively,bauschke2016descent}.
\begin{lemma}\cite[Lemma~1]{bauschke2016descent}\label{lemma:nolip}
	The pair of functions $(g,h)$ is relatively smooth if and only if for all $\x$ and $\y$, it holds that
	\begin{align}
	g(\y)- g(\x) - \inn{\nabla g(\x)}{\y-\x} \leq LD_h(\y,\x).\label{lemma:nolip eq0}
	\end{align}
\end{lemma}
\begin{remark}\label{remark:nolip}
	When $h=\frac{1}{2}\norm{\cdot}^2$, the classical descent lemma is recovered, \ie,
	$
	g(\y)- g(\x) - \inn{\nabla g(\x)}{\y-\x} \leq \frac{L}{2}\norm{\y-\x}^2.\nonumber
	$
\end{remark}

To use Lemma~\ref{lemma:nolip}, we additionally make the following assumptions for the rest of this paper.
\begin{assump}\label{assume: convex}
	The functions $(f_i, h_i)$ are relatively smooth with constants $L_i > 0,\forall i$.
\end{assump}
With the relative smoothness between $(f_i, h_i)$, the following result shows the basic descent property of the proposed method.
\begin{lemma}\label{lemma:descent}
	For any $\x$, and any $i\in \{1,2, \cdots, n\}$, let $\x^+$ to be defined as in E.q. \eqref{eq:x^+}. Then we have
	\begin{align}
	F(\x^+)\leq F(\x) - \left(\frac{1+\theta_i}{\alpha} - L_i\right) D_h(T_i(\x), \x_i),\label{lemma:block descent}
	\end{align}
	where $\theta_i = \theta(h_i)$. In particular, with $ 0 < \alpha < \frac{1+\theta_i}{L_i}$, a sufficient descent in the objective value of $F$ is guaranteed.
\end{lemma}

Maximizing the function $g(\alpha) = (1+\theta_i - L_i\alpha)\alpha$ with respect to $\alpha$ yields the stepsize $\alpha^* = \frac{1+\theta_i}{2L_i}$. Substituting the obtained stepsize into \eqref{lemma:block descent} yields the following result.
\begin{corollary}\label{cor:descent}
	For any $\x$, let $\x^+$ to be defined as in E.q. \eqref{eq:x^+}. With stepsize $\alpha = \frac{1+\theta_i}{2L_i}$, we have
	\begin{align}
	\vspace{-10px}
	F(\x^+)\leq F(\x) - L_i D_h(T_i(\x), \x_i).\label{eq:descent}
	\end{align}
\end{corollary}

With the stepsize $\alpha = \frac{1+\theta_i}{2L_i}$, Corollary~\ref{cor:descent} quantifies the descent in the objective value. Therefore, the stepsize $\alpha^k = \frac{1 + \theta_{i_k}}{2L_{i_k}}$ is an appropriate choice for Algorithm~\ref{alg}.

Since only one block is selected and updated per iteration, the quantity $D_h(\x^+,\x)$ introduced in \cite{lu2018relatively,bauschke2016descent} cannot be used to measure the optimality of the RBCD method. Given an estimate $\x$,  we introduce the reference function $H$ and the corresponding Bregman mapping as follows:
\begin{align}
&H(\x) = \sum_{i=1}^n L_i h_i(\x_i),\\
	&D_H(\y, \x) =\sum_{i=1}^n L_i D_h (\y_i, \x_i)\\
&T(\x) = \argmin_{\u} \inn{\nabla f(\x)}{\u-\x} + D_H(\u, \x) + r(\u). \label{eq:Tx}
\end{align}
Based on this mapping, the following result shows that the quantity $D_H(T(\x), \x)$ can be used to measure the optimality of $F$.

\begin{lemma}\label{lemma:stationary}
	A vector $\x$ is a stationary point of $F$ if and only if $D_H(T(\x), \x) = 0$.
\end{lemma}
Clearly, when $F$ is convex, then the current estimate $\x$ is a global minimum if $D_H(T(\x) , \x) = 0$.

\subsection{Convex and strongly convex case}
In this subsection, we provide the convergence analysis for the case where $F$ is convex. Since $r$ is convex, we have $f$ is also convex. We use $\E_{i}$ (or $\E_{i_k}$) to denote the expectation with respect to a single random variable $i$ (or $i_k$). We use $\E$ to denote the expectation with respect to all random variables $\{i_0, i_1, \cdots\}$.

Instead of using the classical convexity definition, we here use the \textit{relative strongly convexity} introduced in \cite{lu2018relatively}, which is similar to the relative smoothness.
\begin{define}[Relative Strongly Convexity]\cite[Definition 1.2.]{lu2018relatively}
	A function $g$ is $\mu$-strongly convex relative to $h$ if for any $\x$ and $\y$, there exists a scalar $\mu\geq 0$ such that
	\begin{align}
	g(\y)\geq g(\x) + \inn{\nabla g(\x)}{\y-\x} + \mu D_h(\y, \x).
	\end{align}
\end{define}
Note that if $\mu = 0$, the classical convexity for a smooth function $g$ is recovered. Moreover, when $h=\frac{1}{n}\norm{\cdot}$, the classical strongly convexity is recovered. In the rest of this subsection, we assume $f$ is strongly convex relative to $H$.
\begin{assump}
	$f$ is $\mu$-strongly convex relative to $H$, \ie, there exists a scalar $\mu\geq 0$ such that for every $\y$ and $\x$
	\begin{align}
	f(\y)\geq f(\x) + \inn{\nabla f(\x)}{\y-\x} + \mu D_H(\y, \x).\label{eq:relative strongly convex}
	\end{align}
\end{assump}
Since $r$ is assumed to be convex, the function $F$ is also $\mu$-strongly convex relative to $H$, \ie,
\begin{align}
F(\y)\geq F(\x) + \inn{\v}{\y - \x} + \mu D_H(\y, \x),
\end{align}
for some $\v\in \partial F(\x)$. Moreover, by Assumption~\ref{assume: convex}, we have
\begin{align}
	f(T_i(\x) )\leq f(\x) + \inn{\nabla_i f(\x)}{T_i(\x)-\x_i} + L_i D_h(T_i(\x), \x_i).\label{eq:relative smooth}
\end{align}
Substituting $\y = T_i(\x)$ in E.q. \eqref{eq:relative strongly convex} and combing it with the inequality \eqref{eq:relative smooth}, we immediately obtain that $\mu\leq 1$.

The following lemma provides the key inequalities used to prove the convergence results of the RBCD method.
\begin{lemma}\label{lemma:key inequalities}
	For any vector $\x$, let $\x^+$ to be defined as in E.q. \eqref{eq:x^+} by picking up $i\in \{1,2,\cdots, n\}$ uniformly at random. Set stepsize $\alpha = \frac{1+\theta_i}{2L_i}$. For any vector $\u$, the expectation of $F(\x^+)$ satisfies
	\begin{align}
		\E_i[F(\x^+)] \leq& \frac{1}{n} \Big[(n-1) F(\x) + F(\u) \nonumber \\
		&+ (1-\mu) D_H(\u,\x)- D_H(\u, T(\x)) \Big],\label{eq:lemma EF}
	\end{align}
	and the expectation of $D_H(\x^+,\x)$ satisfies
	\begin{align}
	\E_i [D_H(\u, \x^+)] = \frac{n-1}{n} D_H(\u,\x) + \frac{1}{n} D_H(\u, T(\x))\label{eq:lemma ED}.
	\end{align}
\end{lemma}
By applying Lemma~\ref{lemma:key inequalities}, the main convergence results are established in Theorem~\ref{thm:convex}. Note that this result generalizes  \cite[Theorem 1]{lu2015complexity} through replacing the proximal mapping by the Bregman proximal mapping so that the assumption of global Lipschitz-continues (partial) gradient is not necessary.
\begin{theorem}\label{thm:convex}
	Let $\{\x^k\}$ be the sequence generated by Algorithm~\ref{alg}. Then for any $k\geq 0$, the iterates $\x^k$ satisfies
	\begin{align}
		&\E[F(\x^{k}) - F(\x^*)]  \nonumber\\
		\leq& \frac{n}{n+k}\Big(F(\x^*) - F(\x^0) + D_H(\x^*,\x^0)\Big).
	\end{align}
	Further, if $f$ is $\mu$-strongly convex relative to $H$, then
	\begin{align}
		&\E[F(\x^k) - F(\x^*)] \nonumber \\
		\leq &
		\left(1-\frac{(1+\theta)\mu}{n(1+\theta\mu) } \right)^k\Big(F(\x^0)- F(\x^*) + D_H(\x^*,\x^0)\Big),
	\end{align}
	where $\theta = \underset{i}{\min} \{\theta_i\}$.
\end{theorem}
Therefore, if $F$ is convex, the sequence $\{\x^k\}$ needs at most $\bigo{n\varepsilon^{-1}}$ to converge to an $\varepsilon$-solution. Further, the classical linear convergence rate is obtained if $f$ is strongly convex (relative to $H$).

\subsection{Nonconvex case}
In this subsection, we establish the convergence results for the case where $F$ is nonconvex. Since $r$ is convex, $f$ is nonconvex. Due to the nonconvexity, it is of interest to find a stationary point. Lemma~\ref{lemma:stationary} implies that $D_H(T(\x),\x)$ can be used to measure the optimality. The following result shows the descent property of the proposed method in terms of the optimality gap $D_H(T(\x),\x)$.

\begin{lemma}\label{lemma:optimality}
	For any $\x$, let $\x^+$ to be defined as in E.q.\eqref{eq:x^+} by picking up the index $i$ uniformly at random. Let $\alpha = \frac{1+\theta_i}{2L_i}$. Then the following inequality holds:
	\begin{align}
	\E_{i} [F(\x^{+})]\leq F(\x) - \frac{1}{n}D_H(T(\x), \x).\label{eq:nonconvex descent}
	\end{align}
\end{lemma}

Using Lemma~\ref{lemma:optimality}, we can establish the convergence results of the RBCD method for nonconvex $F$.
\begin{theorem}\label{thm:nonconvex}
	Let $\{\x^k\}$ to be the sequence generated by Algorithm~\ref{alg}. Let stepsize $\alpha^k = \frac{1+\theta_{i_k}}{2L_{i_k}}$, then
	\begin{itemize}
		\item[(i)] The sequence $\{F(\x^k)\}$ is non-increasing.
		\item[(ii)] $\sum_{l=0}^\infty \E [D_H(T(\x^l),\x^l) ] <\infty$, and hence the sequence $\{\E [D_H(T(\x^l) ,\x^l)] \}$ converges to zero.
		\item[(iii)]
		$\forall k\geq 0$, we obtain
		\vspace{-5px}
		\begin{align}
			\min_{0\leq l \leq k} \E\left[D_H(T(\x^l), \x^l)\right]
			\leq \frac{n}{k+1}(F(\x^0) - F^*),
		\end{align}
		where $F^*=\inf F(\x) > -\infty$.
		\item[(iv)] Every limit point of $\{\x^k\}$ is a stationary point.
	\end{itemize}
\end{theorem}
Suppose $H$ is $\sigma$-strongly convex with respect to the Euclidean norm $\norm{\cdot}$. Then we have $D_H(\y,\x)\geq\frac{\sigma}{2}\norm{\y-\x}^2$. Combining the strongly convexity of $H$ with Theorem~\ref{thm:nonconvex}, we immediately obtain the following convergence rate result
\begin{align}
\min_{0\leq l\leq k} \E\norm{T(\x^l)- \x^l}^2\leq \frac{2n}{\sigma (k+1)} (F(\x^0) -F^*).
\end{align}
Therefore, the sequence $\{\x^k\}$ converges to a stationary point at the rate of $\bigo{\frac{\sqrt{n}}{\sqrt{k}}}$. In another word, to obtain an $\varepsilon$-stationary point, \ie, $\norm{T(\x) - \x}\leq \varepsilon$, the RBCD method needs to run $\bigo{n\varepsilon^{-2}}$ iterations.

\section{Accelerated Randomized Bregman Coordinate Descent}
In this section, we restrict ourselves to the unconstrained smooth minimization problem as follows
\begin{align}
\minimize_{\x\in\cX}\; f(\x),
\end{align}
where $f$ is convex and satisfies Assumption~\ref{assume:basic}. The closed convex set $\cX$ satisfies $\cX = \cX_1\times \cdots \times \cX_n$ such that $\x_i\in \cX_i$ $\forall i$. It is equivalent to consider $r_i$ as an indicator function of the closed convex set $\cX_i$.

The accelerated randomized Bregman coordinate descent (ARBCD) method is given as Algorithm~\ref{alg:arbcd}. At the $k$-th iteration, the ARBCD method selects a coordinate $i_k$ uniformly at random, and generates the three vectors $\y^k$, $\z^{k+1}$, and $\x^{k+1}$, where the vectors $\y^k$ and $\x^{k+1}$ are the affine combinations of $\x^k$ and $\z^k$, and $\y^k$, $\z^{k}$, and $\z^{k+1}$, respectively, and the vector $\z^{k+1}$ is obtained as follows
\begin{align}
\z^{k+1} = \argmin_{\u\in \cX} \inn{\nabla_{i_k} f(\y^k)}{\u_{i_k} - \y_{i_k}^k} + (n\beta_k)^{\gamma-1}D_H(\u,\z^k). \label{eq:z^+}
\end{align}
Note that Step 1 and 3 of Algorithm~\ref{alg:arbcd} need $\bigo{N}$ operations, while $\bigo{1}$ operations are usually expected in a general coordinate descent method. In the latter section, we will show an efficient implementation of the ARBCD method so that the ARBCD method only needs $\bigo{1}$ operations at each iteration.

\begin{algorithm}[h]
	\caption{Accelerated Randomized Bregman (Block) Coordinate Descent (ARBCD).}
	\label{alg:arbcd}
	\textbf{Input:} initial $\x_0$ and $\gamma$\\
	Initialize: $\z^0 = \x^0$ and $\beta_0 = 1$\\
	\For{$k=1,2\cdots$}{
		\begin{enumerate}[leftmargin=0cm]
			\item $\y^k = (1-\beta_k)\x^k + \beta_k \z^k$
			\item Choose $i_k\in \{1,2,\cdots, n\}$ uniformaly at random \\
			Compute $\z^{k+1}$ by E.q. \eqref{eq:z^+}
			\item $\x^{k+1} = \y^k + n\beta_k (\z^{k+1} - \z^k)$
			\item Choose $\beta_{k+1}\in (0,1]$ such that $\frac{1-\beta_{k+1}}{\beta_{k+1}^\gamma}\leq \frac{1}{\beta_k^\gamma}$
		\end{enumerate}
	}
\end{algorithm}

\section{Convergence Analysis of ARBCD}
To better understand the proposed method, we make the following definitions and observations. First, we define the vector $\tilde{\z}^{k+1}$ as follows
\begin{align}
\tilde{\z}^{k+1} = \argmin_{u\in \cX} \inn{\nabla f(\y^k)}{\u - \y^k} + (n\beta_k)^{\gamma-1}D_H(\u,\z^k), \label{eq:zhat^+}
\end{align}
which is the full-dimensional update version of $\z^{k+1}_{i_k}$ in E.q. \eqref{eq:z^+}. Therefore, the vector $\z^{k+1}$ can be computed by
\begin{align}\label{eq:z^++}
\z_i^{k+1} = \begin{cases}
\tilde{\z}_i^{k+1}, &\text{if $i=i_k$},\\
\z_i^k, &\text{if $i\neq i_k$}.
\end{cases}
\end{align}
It follows from the definition of $\x^{k+1}$ in Step 3 of Algorithm~\ref{alg:arbcd} that we have
\begin{align}\label{eq:x^+ and y}
\x^{k+1}_i = \begin{cases}
\y_i^k + n\beta_k (\z_i^{k+1} - \z_i^k), &\text{if $i=i_k$},\\
\y_i^k, &\text{if $i\neq i_k$}.
\end{cases}
\end{align}
Clearly, the vector $\x^{k+1}$ and $\y^k$ are only one coordinate part from each other, which satisfies the relative smoothness property in Assumption~\ref{assume: convex}.

One of the challenges to establish the convergence results is from the nature of Bregman distances. Since a Bregman distance is in general not a norm, it does not hold the \textit{homogeneous translation invariant}, \ie,
\begin{align}
\norm{\u + \theta(\v-\w)} = \abs{\theta}\norm{\v-\w}, \quad\forall \alpha, \u,\v,\w.
\end{align}
To handle this issue, \cite{hanzely2018accelerated} introduces the notion of \textit{triangle scaling property} (TSP).
\begin{define}\cite[Definition~2]{hanzely2018accelerated}
	The Bregman distance defined with a convex reference function $h$ has the triangle scaling property if there exists some scalar $\gamma > 0$ such that for all $\u,\v,\w$,
	\begin{align}
	D_h((1-\theta)\u + \theta \v, (1-\theta)\u + \theta \w) \leq \theta^{\gamma} D_h(\v,\w), \forall \theta\in [0,1].
	\end{align}
\end{define}
In contrast, we introduce the more general notion of the \textit{generalized translation invariant} (GTI) in the following definition, and show it is equivalent to triangle scaling property, when restricting $\theta\in [0,1]$.
\begin{define}\label{def:GTI}[Generalized Translation Invariant]
	The Bregman distance defined with a convex reference function $h$ has the generalized translation invariant property if there exists some scalar $\gamma \geq 0$ such that for all $ \u, \v, \w$
	\begin{align}
	D_h(\u + \theta(\v-\w),\u)\leq \abs{\theta}^{\gamma}D_h(\v,\w).\quad\forall \theta \in\reals.
	\end{align}
\end{define}
\begin{lemma}\label{lemma:generalized translation invariant}
	The Bregman distance has the generalized translation invariant with $\theta\in[0,1]$ if and only if it holds the triangle scaling property.
\end{lemma}
\begin{remark}\label{remark:GNI}
	Here we gives three examples to show the existences of GNI in some Bregman divergences, while the proof is included in Appendix.
	\begin{itemize}
		\item[(i)] The norms. Let $\norm{\cdot}_\A$ be a norm, $\A$ be a positive define matrix, $h(\x)=(1/2)\norm{\x}^2_\A$, and $D_h(\x,\y)=(1/2)\norm{\x-\y}_\A^2=(1/2)\x^T\A\y$. It is easy to see that $\gamma=2$.
		\item[(ii)] The Kullback-Leibler (KL) divergence. Let $h$ be the negative  Boltzmann-Shannon entropy: $h(\x)=\sum_{i=1}^N \x_i\log \x_i$ defined over $\reals_{+}^N$. The Bregman distance is given by
		\begin{align}
		D_{\KL}(\x,\y) =\sum_{i=1}^N \left(\x_i \log \left(\frac{\x_i}{\y_i} \right) - \x_i + \y_i\right).
		\end{align}
		It can be shown that $\gamma = 1$.
		\item[(iii)] The Itakura-Saito (IS) distance. Let $h$ be the Burg's entropy: $h(\x) = -\sum_{i=1}^N \log\x_i$ on $\reals_{++}^N$. The Bregman distance associated with $h$ is given by
		\begin{align}
		D_{\IS}(\x,\y) =\sum_{i=1}^M \left(- \log \left(\frac{\x_i}{\y_i} \right) + \frac{\x_i }{\y_i}-1\right).
		\end{align}
		To satisfy the definition of GNI, we must have $\gamma = 0$. Similar to TSP, however, $\gamma=0$ is the uniform value for $D_{\IS}$, and the intrinsic $\gamma$ value can be $2$ if the three points are close to each other \cite[Theorem~1]{hanzely2018accelerated}.
	\end{itemize}
\end{remark}

Note that the GTI is more general since TSP needs $\theta \in [0,1]$, but GTI holds for all $\theta\in\reals$.

To use the notion of GTI, we make the following assumption.
\begin{assump}\label{assump:generalized translation invariant}
	The Bregman distances $D_h(\cdot,\cdot)$ have the generalized translation invariant with the constant $\gamma > 0$, $\forall i$.
\end{assump}
Using the notion of GTI, we will show that the ARBCD method converges with a sublinear rate of $\bigo{n\varepsilon^{-1/\gamma}}$. We start with recalling the critical lemma \cite[Lemma~3.2]{chen1993convergence} for a Bregman proximal mapping.
\begin{lemma}\cite[Lemma~3.2]{chen1993convergence}\label{lemma:prox}
	For a convex function $\phi$ and a vector $x$, if the Bregman proximal mapping is defined as
	\begin{align}
	\x^+=\argmin\; \phi(\u) + D_h(\u, \x),
	\end{align}
	and then
	\begin{align}
	\phi(\u) + D_h(\u, \x)\geq \phi(\x^+) + D_h(\x^+, \x) + D_h(\u, \x^+), \forall \u.
	\end{align}
\end{lemma}
The key relationship between two consecutive iterates in Algorithm~\ref{alg:arbcd} is established in the following lemma.
\begin{lemma}\label{lemma:acc ineq}
	Suppose Assumptions~\ref{assume:basic}, \ref{assume: convex}, and \ref{assump:generalized translation invariant} holds. For any vector $\u$, the sequences generated by Algorithm~\ref{alg:arbcd} satisfy, for all $k\geq 0$,
	\begin{align}
	&\E_{i_k} \left[\frac{1-\beta_{k+1}}{\beta_{k+1}^{\gamma}}(f(\x^{k+1}) - f(\u)) + n^\gamma D_H(\u, \z^{k+1}) \right] \nonumber\\
	\leq& \frac{1-\beta_k}{\beta_k^{\gamma}}( f(\x^k) - f(\u)) + n^\gamma D_H(\u,\z^k).
	\end{align}
\end{lemma}
The following lemma introduces a sequence $\{\beta_k\}$ that satisfies the condition in Step 4 of Algorithm~\ref{alg:arbcd}.
\begin{lemma}\cite[Lemma~3]{hanzely2018accelerated}\label{lemma:theta}
	The sequence $\beta_k = \frac{\gamma}{k+\gamma}$ satisfies
	\begin{align}
	\frac{\beta_{k+1} - 1}{\beta_{k+1}^\gamma}\leq \frac{1}{\beta_k^\gamma},\quad\forall k\geq 0.
	\end{align}
\end{lemma}
Combing Lemma~\ref{lemma:acc ineq} with Lemma~\ref{lemma:theta}, the main convergence results for the ARBCD are established in the following theorem.
\begin{theorem}\label{thm:acc}
	Suppose Assumptions~\ref{assume:basic}, \ref{assume: convex}, and \ref{assump:generalized translation invariant} hold. If $\beta_{k} = \frac{\gamma}{k+\gamma}$ for all $k\geq0$, then the following inequality holds, for any vector $\u$,
	\begin{align}
	\E \left[f(\x^{k+1}) - f(\u) \right] \leq \left(\frac{n\gamma}{k+\gamma}\right)^\gamma D_H(\u,\x^0), \;\forall k\geq0.
	\end{align}	
\end{theorem}
Note that due to the affine combinations in Step 1 and 3 of Algorithm~\ref{alg:arbcd}, the current implementation requires $\bigo{N}$ operations. In the next section, we introduce an efficient implementation so that only $\bigo{1}$ operations are needed at each iteration.

\section{Efficient implementation}
In order to avoid full-dimensional vector operations, the previous works \cite{lee2013efficient,fercoq2015accelerated} propose a strategy that changes the variables for the accelerated coordinated descent methods in the global Lipschitz-continuous (partial) gradient setting. Here we show this scheme can be adapted so that the full-dimensional operations can be avoided in the relative smoothness setting, which is given as Algorithm~\ref{alg:efficient}. Instead of computing the vector $\z^{k+1}$, a search direction $\d^k_{i_k}$ is computed in Algorithm~\ref{alg:efficient} as follows
\begin{align}
\d_{i_k}^{k} = &\argmin_{\v_{i_k}^k+\d\in \cX_{i_k}} \inn{\nabla_{i_k} f(\beta_k^\gamma \u^k + \v^k )}{\d} \nonumber\\
&+ (n\beta_k)^{\gamma-1}L_{i_k}D_h(\v^k_{i_k}+\d,\v_{i_k}^k). \label{eq:d^k}
\end{align}
\begin{algorithm}[h]
	\caption{Efficient implementation of ARBCD.}
	\label{alg:efficient}
	\textbf{Input:} initial $\x_0$ and $\gamma$\\
	Initialize: $\v^0 = \x^0$, $\u^0 = \textbf{0}$ and $\beta_0 = 1$\\
	\For{$k=1,2\cdots$}{
		\begin{enumerate}[leftmargin=0cm]
			\item Choose $i_k\in \{1,2,\cdots, n\}$ uniformaly at random \\
			Compute $\d_{i_k}^{k}$ by Eq.\eqref{eq:d^k}
			\item $\v_{i_k}^{k+1} = \v_{i_k}^k + \d^k_{i_k}$
			\item $\u_{i_k}^{k+1} = \u_{i_k}^k -\frac{1-n\beta_{k}}{\beta_{k}^\gamma} \d^k_{i_k}$
			\item Compute $\beta_{k+1}$ from $\frac{1-\beta_{k+1}}{\beta_{k+1}^\gamma}= \frac{1}{\beta_k^\gamma}$
		\end{enumerate}
	}
	\Return{$\beta_{k+1}u_{k+1} + v_{k+1}$}
\end{algorithm}

\begin{prop}\label{prop:eff}
	The sequences $\{\x^k,\y^k,\z^k\}$ and $\{\u^k, \v^k\}$ generated from Algorithm~\ref{alg:arbcd} and \ref{alg:efficient}, respectively, satisfy
	\begin{align}
	\z^k = & \v^k\\
	\x^k = &\beta_{k-1}^\gamma\u^k + \v^k\\
	\y^k = &\beta_k^\gamma \u^k + \v^k,
	\end{align}
	for all $k\geq 1$. That is, these two algorithms are equivalent.
\end{prop}
Note that in Algorithm~\ref{alg:efficient}, only a single block coordinates of the vectors $\u^k$ and $\v^k$ are updated at each iteration, which cost $\bigo{N_i}$ operations. Although computing the partial gradient in E.q. \eqref{eq:d^k} may still cost full-dimensional operations in general, the previous works \cite{lee2013efficient,fercoq2015accelerated,lin2014accelerated} introduce a number of optimization problems where the partial gradient can be computed cheaply without actually forming $\y^k$.

\section{Numerical Experiments}
To showcase the strength of the proposed methods, we consider two applications of relatively smooth convex optimization: Poisson inverse problem, and relative-entropy nonnegative regression.

\subsection{Poisson linear inverse problem}\label{sec:PLI}
A large number of problems in nuclear medicine, night vision, astronomy and hyperspectral imaging can be described as inverse problems where data measurements are collected according to a Poisson process whose underling intensity function is indirectly related to an object of interest through a linear system. This class of problems have been studied intensively in the literature. See for instance \cite{csiszar1991least,dempster1977maximum,engl1996regularization} and references therein, as well as a more recent comprehensive review \cite{bertero2009image} for the up-to-date references.

Formally, in a Poisson inversion problem we are given a nonnegative observation matrix $\A\in\reals_+^{M\times N}$, a noisy measurement vector $\b\in\reals_+^{M}$, and the goal is to recover the signal or image of interest $\x\in \reals_+^N$. Under the Poisson assumption, we can rewrite the observation model as follows
\begin{align}
	\b \sim\text{Poisson}(\A\x).
\end{align}
Therefore, a natural and widely used measure of proximity of two nonnegative vectors is based on the KL divergence. Particularly, minimizing the KL-divergence $D_\KL(\b, \A\x)$ is equivalent to maximize the Poisson log-likelihood function. The optimization problem can be formulated as follows
\begin{align}
	\minimize_{\x\geq 0} \; f(\x) \equiv D_\KL(\b, \A\x).
\end{align}
To apply the RBCD and ARBCD methods, we need to identify a series of adequate reference functions $h_i$. Here we use Burg’s entropy and the corresponding Bregman distance, \ie, the IS distance. 
\begin{lemma}\label{lemma:relative smooth}
	Let $f_i(\x_i) = D_\KL(\b, \A\x)$ and $h_i(\x_i)$ to be defined as
	\begin{align}
		h_i(\x_i) = -\log\x_i.
	\end{align}
	Then the functions $(f_i, h_i)$ are relatively smooth with any scalar $L_i$ satisfying
	\begin{align}
		L_i\geq \norm{\b}_1=\sum_{i=1}^M \b_i.
	\end{align}
\end{lemma}
Equipped with Lemma~\ref{lemma:relative smooth}, Theorem~\ref{thm:convex} is applicable and warrants the convergence. Since $\theta(h_i)=0$, we can take the stepsize $\alpha_k = \frac{1}{2\norm{\b}_1}$, $\forall k\geq0$. To solve Poisson inverse problems, the E.q. \eqref{eq:T_i} can be written as
\begin{align}
	T_i(\x) = \argmin_{\u_i\geq 0} \inn{\nabla f_i(\x_i)}{\u_i} + 2\norm{\b}_1 D_\IS(\u_i, \x_i).
\end{align}
It follows from \cite[Theorem~1]{hanzely2018accelerated} that the intrinsic TSE of a Bregman distance is $2$, even the uniform TSE is not. In addition, \cite{hanzely2018accelerated} numerically shows the convergence and efficiency of the Accelerated Bregman Proximal method (ABPG) with $\gamma=2$. Thus, we here also use $\gamma=2$ for the ARBCD method. As a result, E.q. \eqref{eq:z^+} becomes
\begin{align}
\z_{i_k}^{k+1} = \argmin_{\u_{i_k}\geq 0} \inn{\nabla_{i_k} f(\y^k)}{\u_{i_k}} + 2n\beta_k \norm{\b}_1 D_\IS(\u_{i_k}, \z^k_{i_k}).
\end{align}

We compare the proposed algorithms RBCD and ARBCD with two state-of-the-art algorithms: Bregman Proximal Gradient (BPG) method \cite{bauschke2016descent} and accelerated Bregman Proximal Gradient (ABPG) \cite{hanzely2018accelerated} method. All algorithms are implemented in Matlab code.

\begin{figure}[t]
	\centering
		\centerline{\includegraphics[width=0.8\columnwidth]{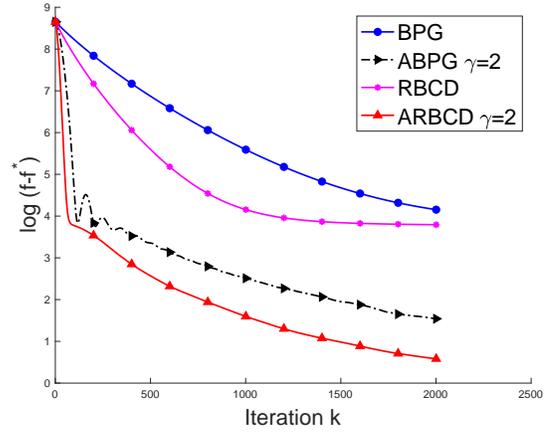}}
	\caption{Poisson inverse problem: synthetic dataset with $M=500$ and $N=500$.}
			\label{fig:iterations}
\end{figure}

\begin{figure}[t]
	\centering
	\centerline{\includegraphics[width=0.8\columnwidth]{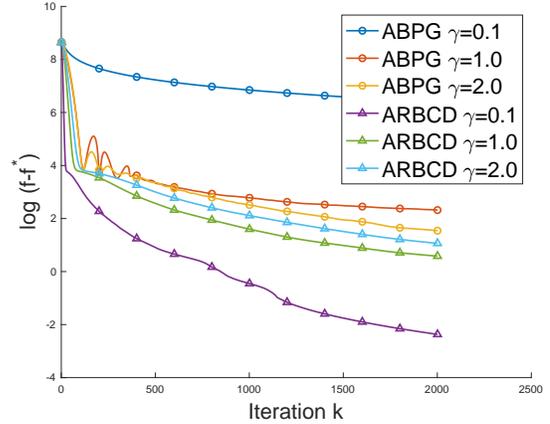}}
	\caption{Poisson inverse problem: synthetic dataset with varying $\gamma$ values.}
	\label{fig:gamma}
\end{figure}

Figure~\ref{fig:iterations} shows the computational results for a randomly generated dataset with $M=500$ and $N=500$. The entries in $\A$ and $\b$ are generated randomly from a uniform distribution over the interval $[0,1]$. Each algorithm starts with the same initial values. Note that the CD-type methods has a inner loop of $N$ iterations as their computational complexity is $N$ times cheaper than the gradient-based methods. As a result, the computational complexity in each iteration is identical.

In Figure~\ref{fig:iterations}, we can see the RBCD method is only slightly better than the BPG method, because the RBCD method uses the most updated coordinate to update, and BPG and RBCD methods use the same stepsize $\alpha_k= \frac{1}{2\norm{b}_1}$. Figure~\ref{fig:iterations} also shows that the accelerated methods ABPG and ARBCD are both faster than their non-accelerated variants. We can also conclude that the ARBCD method is faster than the other methods. It is well-known that the accelerated (proximal) gradient method does not guarantee the descent in the objective values at each iteration. Instead, the number of ripples are on the traces of the objective values. This criteria can be found on the ABPG method as well in Figure~\ref{fig:iterations}. On the other hand, we does not find such ripples or bumps from the ARBCD method. Particularly, Figure~\ref{fig:iterations} shows that the ARBCD method provides consistent descent in the objective values.

It is easy to check numerically that $D_\IS$ does not hold GTI or TSP property for any scalar $\gamma > 0.5$. We conduct another experiment to explore the impact of the parameter $\gamma$. Figure~\ref{fig:gamma} shows the convergence behaviors of the ABPG and ARBCD methods with $\gamma = 0.1, 1.0$ and $2.0$. The larger $\gamma$ is, the more acceleration the ABPG method obtains. However, it seems the ARBCD method holds the opposite relationship with the $\gamma$ values. The ARBCD method achieves the maximum acceleration when the $\gamma$ is minimum.

\subsection{Relative-entropy nonnegative regression}
Anther formulation to solve the nonnegative linear inverse problem introduced in Section~\ref{sec:PLI} is to minimize $D_{\KL}(\A\x,\b)$, \ie,
\begin{align}
	\minimize_{\x\geq 0}\; f(\x)\equiv D_\KL(\A\x, \b).
\end{align}
In this case, the following result shows that the function $f$ is relative smooth to the Boltzman-Shannon entropy defined by
\begin{align}
	h(\x) = \x\log \x, \quad\forall \x\in\reals_{+}.
\end{align}
\begin{lemma}\label{lemma:relative entropy}
	Let $f_i(\x_i) = D_\KL(\A\x,\b)$ and $h_i(\x_i)$ to be defined as
	\begin{align}
	h_i(\x_i) = \x_i\log\x_i.
	\end{align}
	Then the functions $(f_i, h_i)$ are relatively smooth with any scalar $L_i$ satisfying
	\begin{align}
	L_i\geq \sum_{i=1}^M \a_{ij},
	\end{align}
	where $\a_{ij}$ is the $(i,j)$-th entry of $\A$.
\end{lemma}

\begin{figure}[t]
	\centering
	\centerline{\includegraphics[width=0.8\columnwidth]{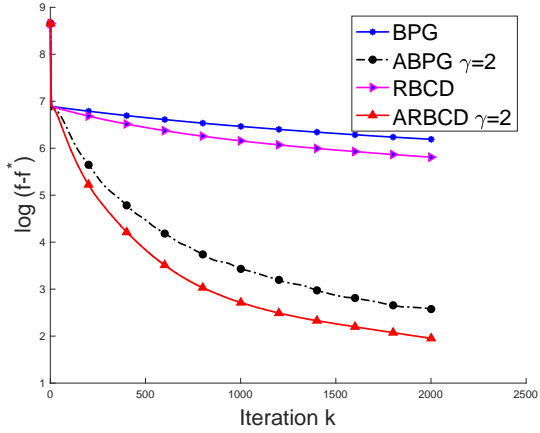}}
	\caption{Relative-entropy nonnegative regression: synthetic dataset with $M=500$ and $N=500$.}
	\label{fig:iterations_RE}
\end{figure}

\begin{figure}[t]
	\centering
	\centerline{\includegraphics[width=0.8\columnwidth]{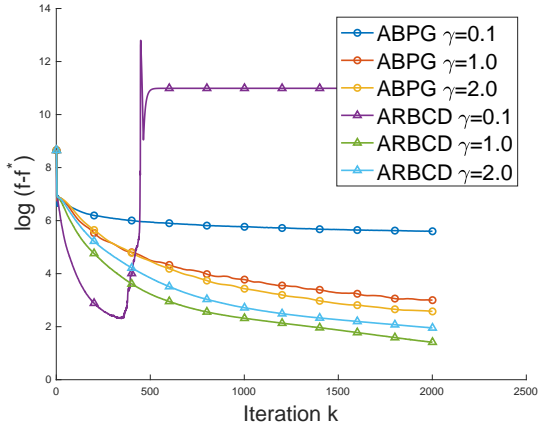}}
	\caption{Relative-entropy nonnegative regression: synthetic dataset with varying $\gamma$ values.}
	\label{fig:gamma_RE}
\end{figure}

Figure~\ref{fig:iterations_RE}-\ref{fig:gamma_RE} shows the computational results for a randomly generated dataset with $M=500$ and $N=500$. Figure~\ref{fig:iterations_RE} shows the almost identical convergence behaviors as in Figure~\ref{fig:iterations}, where the RBCD and ARBCD methods are slightly faster than the BPG and ABPG methods, respectively, and the ARBCD method is faster than the rest methods. As the $\gamma$ values increases, Figure~\ref{fig:gamma_RE} shows improved convergence for the ABPG method. However, the smallest value of $\gamma$, \ie, $\gamma = 0.1$, causes the divergence of the ARBCD method. Therefore, the choice of the hyperparameter $\gamma$ has significant influence on the performance of the ARBCD method.

\section{Conclusion}
In this paper, we propose a randomized Bregman (block) coordinate descent (RBCD) method and its accelerated variant ARBCD method for minimizing a composite problems, where the smooth part of the objective function does not satisfies the global Lipschitz-continuous (partial) gradient property. By using the relative smoothness, we establish the iteration complexity of $\bigo{n\varepsilon^{-2}}$ to obtain an $\varepsilon$-stationary point in the case where $F$ is nonconvex. Besides, the iteration complexity is improved to $\bigo{ n\varepsilon^{-1}}$ if $f$ is convex, and the global linear convergence rate can be achieved by RBCD if $f$ is strongly convex. We introduce the notion of generalized translation invariant. Thanks to this notion, we are able to establish the convergence result for the ARBCD method which uses the acceleration technique. Thus, the iteration complexity is further improved to $\bigo{ n\varepsilon^{-1/\gamma}}$ by the ARBCD method.

\appendix
\section{Appendix}
\subsection{Proof of Lemma~\ref{lemma:descent}}
\begin{proof}
	From the relative smoothness, we obtain
	\begin{align}
	f(\x^{+})\leq f(\x) + \inn{\nabla_i f(\x)}{T_i(\x ) - \x_i} + L_i D_h(T_i(\x), \x_i).\label{eq:descent eq}
	\end{align}
	From the optimality of $T_i(\x)$ in \eqref{eq:T_i}, we have
	\begin{align*}
	\nabla_i f(\x) + \frac{1}{\alpha}\left(\nabla h_i (T_i(\x) - \nabla h_i (\x_i) \right) + \v_i^{+} = 0,
	\end{align*}
	for some $\v_i^{+}\in\partial r_i(T_i(\x) )$. The convexity of $r_i$ implies
	\begin{align}
	&r_i(\x_i)-r_i(T_i(\x) )\geq  \inn{\v_i^{+}}{\x_i - T_i(\x )} \nonumber\\
	=&  - \inn{\nabla_i f(\x) + \frac{1}{\alpha}\left(\nabla h_i (T_i(\x) - \nabla h_i (\x_i) \right)  }{\x_i - T_i(\x )}\nonumber\\
	=&-\inn{\nabla_i f(\x)}{\x_i - T_i(\x )} + \frac{1}{\alpha}\left(D_h(\x_i, T_i(\x)) +D_h(T_i(\x), \x_i) \right)
	\label{eq:convex r_i}
	\end{align}
	Combining E.q. \eqref{eq:descent eq} and \eqref{eq:convex r_i} yields
	\begin{align*}
	&f(\x^{+}) + r_i(T_i(\x) ) \\
	\leq& f(\x) + r_i(\x_i) + L_iD_h(T_i(\x),\x_i ) \nonumber\\
	&- \frac{1}{\alpha}\left(D_h(\x_i, T_i(\x)) +D_h(T_i(\x), \x_i) \right)\\
	\leq& f(\x) + r_i(\x_i) - \left(\frac{1+\theta}{\alpha} - L_i\right) D_h(T_i(\x), \x_i),
	\end{align*}
	where the second inequality is due to $D_h(\x_i, T_i(\x))\geq \theta D_h(T_i(\x), \x_i) $. Since $\x_j^{+}=\x_j$ $\forall i\neq j$, we obtain
	\begin{align*}
	F(\x^{+}) \leq F(\x) - \left(\frac{1+\theta}{\alpha} - L_i\right) D_h(T_i(\x), \x_i).
	\end{align*}
\end{proof}

\subsection{Proof of Lemma~\ref{lemma:stationary}}
\begin{proof}
	$(\Longrightarrow)$. Suppose $\x$ is a stationary point. Then we have
	\begin{align*}
	\nabla f(\x)  + \v = 0,
	\end{align*}
	for some $\v\in \partial r(\x)$. From the convexity of $r$, it follows that for any vector $\u$
	\begin{align}
	r(\u)\geq r(\x) - \inn{\nabla f(\x)}{\u - \x}.\label{eq:opt1st}
	\end{align}
	By the optimality of \eqref{eq:Tx}, we obtain
	\begin{align}
	\nabla f(\x) + \nabla H(T(\x)) - \nabla H(\x) + \v^+ =0,\label{eq:opt3rd}
	\end{align}
	for some $\v^+\in \partial r(T(\x))$. It follows that
	\begin{align}
	r(\x)- r(T(\x))&\geq - \inn{\nabla f(\x)}{\x - T(\x)} \nonumber\\
	&- \inn{\nabla H(T(\x)) - \nabla H(\x)}{\x-T(\x)}.\label{eq:opt2nd}
	\end{align}
	Let $\u = T(\x)$ and combine the equations \eqref{eq:opt1st} and \eqref{eq:opt2nd}. Then we obtain
	\begin{align*}
	0\geq D_H(\x, T(\x)) + D_H(T(\x), \x).
	\end{align*}
	Since $D_H(\x, T(\x)),D_H(T(\x), \x)\geq 0$, we obtain $D_H(T(\x), \x) = 0$.
	
	$(\Longleftarrow)$. Suppose $D_H(T(\x), \x) = 0$. The (strict) convexity of $H$ implies $T(\x) = \x$. From \eqref{eq:opt3rd}, we obtain
	\begin{align*}
	0\in \nabla f(\x) + \partial r(\x),
	\end{align*}
	which indicates $\x$ is a stationary point.
\end{proof}

\subsection{Proof of Lemma~\ref{lemma:key inequalities}}
\begin{proof}
	Since each block $i$ is selected uniformly at random, we have
	\begin{align*}
	&\E_i[F(\x^+)] = \sum_{i=1}^n \frac{1}{n} F(\x^+)\\	
	=&\frac{1}{n}  \sum_{i=1}^n f(\x^+) + r(\x^+)\\
	\overset{(i)}{\leq}&\frac{1}{n}  \sum_{i=1}^n f(\x) + \inn{\nabla_i f(\x)}{T_i(\x) - \x_i} \nonumber\\
	&+ L_i D_h(T_i(\x), \x_i) + r_i(T_i(\x)) + \sum_{j\neq i} r_j(\x_j)\\
	\overset{(ii)}{\leq}&\frac{1}{n}[n f(\x) + \inn{\nabla f(\x)}{T(\x) - \x} \nonumber\\
	&+ D_H(T(\x),\x ) + r(T(\x)) + (n-1) r(\x)]\\
	=&\frac{1}{n} [(n-1) F(\x) + f(\x) + \inn{\nabla f(\x)}{T(\x) - \x}\nonumber\\
	& + D_H(T(\x),\x ) + r(T(\x)) ]\\
	\overset{(iii)}{=}&\frac{1}{n} [(n-1) F(\x) + f(\u) -\mu D_H(\u,\x) \\
	&+ D_H(T(\x), \x)
	+ r(\u) + \inn{\nabla H(T(\x)) - \nabla H(\x) }{\u - T(\x)} ]\\
	\overset{(iv)}{=}&\frac{1}{n} \left[(n-1) F(\x) + F(\u) + (1-\mu) D_H(\u,\x) - D_H(\u, T(\x)) \right]
	\end{align*}
	where $(i)$ follows from the relative smoothness of $(f_i, h_i)$; $(ii)$ uses the fact of $T_i(\x) = T(\x)_i$; $(iii)$ is based on the convexity of $f$ and $r$; $(iv)$ uses the the fact of $\inn{\nabla h(\z) - \nabla h(\x) }{\y - \z}=D_h(\y,\x) -D_h(\y,\z) - D_h(\z,\x)$.
	
	For any vector $\u$, we have
	\begin{align}
	&D_H(\u, \x^+) \nonumber\\
	=& L_i D_h(\u_i, T_i(\x)) + \sum\nolimits_{j\neq i} L_j D_h(\u_j, \x_j)\nonumber\\
	=&L_i D_h(\u_i, T_i(\x)) - L_i D_h(\u_i, \x_i) + D_H(\u, \x) \label{eq:lemma E D_H}
	\end{align}
	Taking the expectation of Eq.\eqref{eq:lemma E D_H} with respect to $i$ yields
	\begin{align*}
	&\E_i [D_H(\u, \x^+)]  \\
	=& \E_i \left[ D_H(\u,\x) - L_i D_h(\u_i, \x_i) + L_i D_h(\u_i, T_i(\x))  \right]\\
	=&\sum_{i=1}^n \frac{1}{n} \left[ D_H(\u,\x) - L_i D_h(\u_i, \x_i) + L_i D_h(\u_i, T_i(\x))  \right]\\
	=& \frac{1}{n} \left[ n D_H(\u,\x) - D_H(\u, \x) + D_H(\u, T(\x))  \right]\\
	=& D_H(\u,\x) - \frac{1}{n} \left[ D_H(\u, \x) - D_H(\u, T(\x))  \right]
	\end{align*}
\end{proof}

\subsection{Proof of Theorem~\ref{thm:convex}}
\begin{proof}
	Combining \eqref{eq:lemma ED} with \eqref{eq:lemma EF}, let $\u=\x^*$, and we have
	\begin{align}
	&\E_i[F(\x^+) + D_H(\x^*, \x^+)] \nonumber\\
	\leq& \frac{n-1}{n} F(\x) + \frac{1}{n} F(\x^*) + \left(1-\frac{\mu}{n} \right) D_H(\x^*, \x)\label{eq:prop descent2}\\
	\leq& \frac{n-1}{n} F(\x) + \frac{1}{n} F(\x^*) + D_H(\x^*, \x).\label{eq:prop descent1}
	\end{align}
	Taking the expectation of \eqref{eq:prop descent1} with respect to $\{i_0, i_1, \cdots\}$ yields
	\begin{align*}
	\E[F(\x^+) ] \leq& \E [F(\x)  + D_H(\x^*, \x) - D_H(\x^*, \x^+)\nonumber\\
	 &- \frac{1}{n}\left(F(\x)  -F(\x^*)\right)].
	\end{align*}
	Summing over $l=0,1,\cdots, k-1$ yields
	\begin{align*}
	&\E[F(\x^{k})] \\
	\leq& F(\x^0) +D_H(\x^*,\x^0)-\E [D_H(\x^*, \x^{k})]\\
	&- \frac{1}{n}\sum_{l=0}^{k-1} \E\left[ F(\x^l) - F(\x^*)\right]\\
	\leq&F(\x^0) +D_H(\x^*,\x^0)- \frac{1}{n}\sum_{l=0}^{k-1} \E\left[ F(\x^l) - F(\x^*)\right]\\
	\leq&F(\x^0) +D_H(\x^*,\x^0)- \frac{k}{n} \E\left[ F(\x^{k+1}) - F(\x^*)\right],
	\end{align*}
	where the last inequality is because $\{F(\x^{l})\}$ is a descent sequence. Subtracting $F(\x^*)$ on both sides and rearrange yields
	\begin{align*}
	\frac{n+k}{n}\E[F(\x^{k}) - F(\x^*)] \leq F(\x^*) - F(\x^0) + D_H(\x^*,\x^0).
	\end{align*}
	Dividing both sides by $\frac{n+k}{n}$ yields the desired result.
	
	If $f$ is $\mu$-strongly convex relative to $H$, we have
	\begin{align*}
	&\E_i[F(\x^+) + D_H(\x^*, \x^+)] \\
	\leq& \frac{n-1}{n} F(\x) + \frac{1}{n} F(\x^*) + \left(1-\frac{\mu}{n} \right) D_H(\x^*, \x).
	\end{align*}
	Subtracting $F(\x^*)$ on the both sides and rearrange yields
	\begin{align}
	&\E_i[F(\x^+) - F(\x^*)+ D_H(\x^*, \x^+)]\nonumber\\
	 \leq& F(\x)- F(\x^*) + D_H(\x^*,\x) \nonumber\\
	 &-\frac{1}{n}F(\x) -  F(\x^*) + \mu D_H(\x^*, \x).\label{eq:prop convex 1st}
	\end{align}
	The relative strongly convexity of $F$ implies
	\begin{align*}
	F(\x) - F(\x^*) + \mu D_H(\x^*,\x)\geq& \mu D_H(\x, \x^*)+ \mu D_H(\x^*,\x)\\
	\geq& (1+\theta) \mu D_H(\x^*,\x).
	\end{align*}
	Define
	\begin{align}
	\beta = \frac{(1+\theta)\mu}{1+\theta \mu}.
	\end{align}
	Clearly, we have $\beta \leq 1$ since $\mu\leq 1$. Then
	\begin{align*}
	&F(\x) - F(\x^*) + \mu D_H(\x^*,\x)\\
	\geq& \beta(F(\x) - F(\x^*) + \mu D_H(\x^*,\x)) + (1-\beta)(1-\theta)\mu D_H(\x^*,\x)\\
	=&\beta(F(\x) - F(\x^*) +D_H(\x^*,\x)).
	\end{align*}
	Combining the inequality above with \eqref{eq:prop convex 1st} yields
	\begin{align*}
	&\E_i[F(\x^+) - F(\x^*)+ D_H(\x^*, \x^+)]\\
	 \leq& \left(1-\frac{\beta}{n} \right)\left(F(\x)- F(\x^*) + D_H(\x^*,\x)\right)
	\end{align*}
	Taking the expectation with respect to $\{i_0, i_1, \cdots\}$ on the both sides of the relation above, we have
	\begin{align*}
	&\E[F(\x^k) - F(\x^*)+ D_H(\x^*, \x^k)]\\
	 \leq& \left(1-\frac{\beta}{n} \right)^k\left(F(\x^0)- F(\x^*) + D_H(\x^*,\x^0)\right).
	\end{align*}
	Dropping $D_H(x^*,\x^k)$ on the left hand yields the desired result.
\end{proof}

\subsection{Proof of Lemma~\ref{lemma:optimality}}
\begin{proof}
	Taking the expectation of \eqref{eq:descent} with respect to $i$ yields
	\begin{align*}
	\E_i [F(\x^{+})]&\leq F(\x) - \E_i[L_i D_h(T_i(\x), \x_i)]\\
	&=F(\x) - \sum_{i=1}^n\frac{1}{n} L_i D_h(T_i(\x), \x_i)\\
	&= F(\x) - \frac{1}{n}\sum_{i=1}^n L_i D_h(T_i(\x), \x_i)\\
	&\overset{(i)}{=} F(\x) - \frac{1}{n}\sum_{i=1}^n L_i D_h(T(\x)_i, \x_i)\\
	&= F(\x) - \frac{1}{n} D_H(T(\x), \x),
	\end{align*}
	where $(i)$ is because $T_i(\x) = T(\x)_i$.
\end{proof}

\subsection{Proof of Theorem~\ref{thm:nonconvex}}
\begin{proof}
	\noindent$(i)$. The result is directly obtained from Lemma~\ref{lemma:optimality}.
	
	\noindent$(ii)$. Taking the expectation of \eqref{eq:nonconvex descent} with respect to all variables and rearranging yields
	\begin{align*}
	\E\left[D_H(T(\x^l) ,\x^l)  \right] \leq n\E \Big(F(\x^l) - F(\x^{l+1}) \Big).
	\end{align*}
	Taking the telescopic sum of the above inequality for $l=0,1,\cdots, k$ gives us
	\begin{align}
	\sum_{l=0}^k \E\left[D_H(T(\x^l) ,\x^l)  \right]
	\leq&  n\left(F(\x^0) - \E [F(\x^{K+1})] \right)\nonumber\\
	\leq& n \left( F(\x^0) - F^* \right).\label{eq:sum}
	\end{align}
	Since $F$ is lower bounded, taking the limit $k\rightarrow\infty$ yields the desired result.
	
	\noindent$(iii)$. The inequality \eqref{eq:sum} further implies that
	\begin{align*}
	(k+1)\min_{0\leq l \leq k} \E\left[D_H(T(\x^l), \x^l)\right]
	\leq& \sum_{l=0}^k \E\left[D_H(T(\x^l) ,\x^l)  \right]\\
	\leq& n(F(\x^0) - F^*).
	\end{align*}
	Dividing $k+1$ on both sides gives us the desired result.
	\noindent $(iv)$. Let $\x^*$ to be a limit point of $\{\x^k\}$ and there exists a subsequence $\{\x^{k_p}\}$ such that $\x^{k_p} \rightarrow \x^*$ as $p\rightarrow \infty$.
	
	Since the functions $r_i$ are lower semi-continuous, we have for all $i$,
	\begin{align}
	\liminf_{p\rightarrow \infty} r_i(\x_i^{k_p}) \geq r_i(\x_i^*).\label{eq:lower semi}
	\end{align}
	At the $k$-th iteration, suppose the index $i$ is selected, then the convexity of $r_i$ implies that
	\begin{align*}
	&r_i(\x_i^{k+1})- r_i(\x_i^*) \\
	\leq&\inn{\nabla_i f(\x^k) + \nabla h_i(\x_i^{k+1}) - \nabla h_i(\x_i^k)}{\x_i^* - \x_i^{k+1}}
	\end{align*}
	Let $\{\x^{k_q}\}$ be the subsequence of $\{\x^{k_p}\}$ such that the index $i$ is selected. Choosing $k = k_q-1$ in the above inequality, and letting $q\rightarrow$ yields
	\begin{align}
	\limsup_{q\rightarrow \infty} r_i(\x_i^{k_q}) \leq r_i(\x^*_i),\label{eq:semi cont}
	\end{align}
	where we use the facts $\x^{k_q}\rightarrow \x^*$ as $q\rightarrow \infty$. Thus, combining \eqref{eq:semi cont} with \eqref{eq:lower semi}, we have
	\begin{align*}
	\lim_{q\rightarrow \infty} r_i(\x_i^{k_q}) = r_i(\x^*_i).
	\end{align*}
	Since $i$ is selected arbitrarily, we have
	\begin{align*}
	\lim_{p\rightarrow \infty} r_i(\x_i^{k_p}) = r_i(\x^*_i), \quad\forall i.
	\end{align*}
	Furthermore, by the continuity of $f$, we obtain
	\begin{align*}
	\lim_{p\rightarrow \infty} F(\x^{k_p})
	=& \lim_{p\rightarrow \infty}\left\{f(\x^{k_p}) + \sum_{i=1}^n r_i(\x^{k_p})\right\}\\
	=&f(\x^*) + \sum_{i=1}^n r_i(\x_i^*) = F(\x^*).
	\end{align*}
	From $(ii)$ and Lemma~\ref{lemma:stationary}, it follows that $\x^*$ is a stationary point of $F$.
\end{proof}

\subsection{Proof of Lemma~\ref{lemma:generalized translation invariant}}
\begin{proof}
	$\Longrightarrow$. Suppose the Bregman distance $D_h(\cdot, \cdot)$ holds the generalized translation variant, and let $\u = (1-\theta)\x + \theta \w$ for any $\x$. Then we have
	\begin{align*}
	D_h((1-\theta)\x +\theta \v, (1-\theta)\x + \theta \w) \leq \abs{\theta}^\gamma D_h(\v,\w),\quad\forall \theta\in\reals.
	\end{align*}
	Since the above inequality holds for all $\theta$, it must hold for $\theta\in [0,1]$.
	
	$\Longleftarrow$. Suppose the triangle scaling property holds. Let $\y = (1-\theta)\u + \theta \w$, then we have
	\begin{align}
	D_h(\y +\theta (\v-\w), \y) \leq \theta^\gamma D_h(\v,\w),\quad\forall \theta\in[0,1].
	\end{align}
	Therefore, the generalized translation invariant holds for $\theta\in [0,1]$.
\end{proof}

\subsection{Proof of Remark~\ref{remark:GNI}}
	\begin{itemize}
		\item[(i)] It is easy to verify that
		\begin{align*}
		\frac{1}{2}\norm{\u + \theta(\v -\w) -\u}^2_\A = \frac{1}{2}\theta^2\norm{\v - \w}_\A^2.
		\end{align*}
		\item[(ii)] Without loss the generality, we assume $N=1$. Using the log sum inequality, we obtain 
		\begin{align*}
		&D_\KL(\u + \theta(\v-\w), \u) \\
		=& (\u + \theta(\v-\w))\log \left(\frac{\u + \theta(\v-\w)}{\u} \right) + \theta(\v-\w) \\
		=& (\u + \theta(\v-\w))\log \left(\frac{\u + \theta(\v-\w)}{\u} \right) + \theta(\v-\w) \\
		&+(\theta \v - \u -\theta(\v-\w)) \log \left(\frac{\theta \v - \u -\theta(\v-\w)}{\theta\w-\u}\right)\\
		&-(\theta \v - \u -\theta(\v-\w)) \log \left(\frac{\theta \v - \u -\theta(\v-\w)}{\theta\w-\u}\right)\\
		\leq& \theta \v\log \left(\frac{\v}{\w}\right) + \theta(\v -\w)\\
		= &\theta D_\KL(\v,\w).
		\end{align*}
		\item[(iii)] Without loss generality, we assume $N=1$. As the GNI property in Definition~\ref{def:GTI} is defined  for all $\u,\v,\w$, we consider a special case of $\u = \theta \w$. Then, we have
		\begin{align*}
			&D_\IS(\u + \theta (\v -\w), \u) = D_\IS(\theta \v, \theta \w) \\
			=&-\log \left(\frac{\theta \v }{\theta \u} + \frac{\theta \v }{\theta \u} - 1\right)\\
			=& D_\IS(\v,\w).
		\end{align*}
		To obtain $D_\IS(\theta \v, \theta \w) \leq \abs{\theta}^\gamma D_\IS(\v,\w)$ for all $\theta\in \reals$, we must have $\gamma = 0$, otherwise $1 > \theta^{\gamma}$ for all $\theta\in (0,1)$.
	\end{itemize}

\subsection{Proof of Lemma~\ref{lemma:acc ineq}}
\begin{proof}
	With simple algebra operations, we have
	\begin{align}
	\x^{k+1} -\y^k= n\left[\beta_k(\z^{k+1} - \y^k) + (1-\beta_k)(\x^k - \y^k)\right].\label{eq:x^+, x}
	\end{align}
	Based on the relation in E.q. \eqref{eq:x^+ and y}, we know $\x^{k+1}$ and $\y^k$ satisfy the relative smoothness property since they are only one coordinate difference from each other. Therefore, we obtain
	\begin{align*}
	f(\x^{k+1}) \leq& f(\y^k) + \inn{\nabla_{i_k} f(\y^k)}{\x_{i_k}^{k+1} - \y^k_{i_k}} + L_{i_k} D_h(\x^{k+1}_{i_k}, \y^k_{i_k})\\
	=& f(\y^k) + \inn{\nabla_{i_k} f(\y^k)}{\x_{i_k}^{k+1} - \y^k_{i_k}} \\
	&+ L_{i_k} D_h(\y_{i_k} + n\beta_k(\z_{i_k}^{k+1}-\z_{i_k}^k), \y^k_{i_k})\\
	\overset{(i)}{\leq}& f(\y^k) + \inn{\nabla_{i_k} f(\y^k)}{\x_{i_k}^{k+1} - \y^k_{i_k}} \\
	&+ (n\beta_k)^{\gamma} L_{i_k} D_h((\z_{i_k}^{k+1},\z_{i_k}^k))\\
	\overset{(ii)}{=}& f(\y^k) + n\beta_k\inn{\nabla_{i_k} f(\y^k)}{\z^{k+1}_{i_k} - \y_{i_k}^k} \\
	&+  n(1-\beta_k)\inn{\nabla_{i_k} f(\y^k)}{\x_{i_k}^k - \y_{i_k}^k}\\
	&+ (n\beta_k)^{\gamma} L_{i_k} D_h((\z_{i_k}^{k+1},\z_{i_k}^k))\\
	\overset{(iii)}{=}&
	\beta_k\left[f(\y^k) + n\inn{\nabla_{i_k} f(\y^k)}{\tilde{\z}^{k+1}_{i_k} - \y_{i_k}^k}  \right]\\
	&+(1-\beta_k)\left[f(\y^k) + n\inn{\nabla_{i_k} f(\y^k)}{\x_{i_k}^k - \y_{i_k}^k} \right]\\
	&+ (n\beta_k)^{\gamma} L_{i_k} D_h((\tilde{\z}_{i_k}^{k+1},\z_{i_k}^k)),
	\end{align*}
	where $(i)$ is using the generalized translation invariant, $(ii)$ is due to E.q. \eqref{eq:x^+, x}, and $(iii)$ is due to E.q. \eqref{eq:z^++}. Taking the expectation with respect to $i_k$ on both sides yields for all $\u$
	\begin{align*}
	&\E_{i_k} f(\x^{k+1}) \\
	\leq&
	\beta_k\left[f(\y^k) + n\E_{i_k}\inn{\nabla_{i_k} f(\y^k)}{\tilde{\z}^{k+1}_{i_k} - \y_{i_k}^k}  \right]\\
	&+(1-\beta_k)\left[f(\y^k) + n\E_{i_k}\inn{\nabla_{i_k} f(\y^k)}{\x_{i_k}^k - \y_{i_k}^k} \right]\\
	&+ (n\beta_k)^{\gamma} \E_{i_k}\left[L_{i_k} D_h((\tilde{\z}_{i_k}^{k+1},\z_{i_k}^k))\right]\\
	\overset{(i)}{=}&
	\beta_k\left[f(\y^k) + \inn{\nabla f(\y^k)}{\tilde{\z}^{k+1} - \y^k}  \right]\\
	&+(1-\beta_k)\left[f(\y^k) + \inn{\nabla f(\y^k)}{\x^k - \y^k} \right]\\
	&+ n^{\gamma-1}\beta_k^{\gamma} D_H(\tilde{\z}^{k+1},\z^k)\\
	\overset{(ii)}{\leq}& (1-\beta_k) f(\x^k)\\
	&+ 	\beta_k[f(\y^k) + \inn{\nabla f(\y^k)}{\tilde{\z}^{k+1} - \y^k} \\
	&+ (n\beta_k)^{\gamma-1} D_H(\tilde{\z}^{k+1},\z^k)]	\\	
	\overset{(iii)}{\leq}& (1-\beta_k) f(\x^k)\\
	&+ 	\beta_k[f(\y^k) + \inn{\nabla f(\y^k)}{\u - \y^k} \\
	&+ (n\beta_k)^{\gamma-1} D_H(\u,\z^k)
	- (n\beta_k)^{\gamma-1} D_H(\u,\tilde{\z}^{k+1})
	]\\
	\overset{(iv)}{\leq}& (1-\beta_k) f(\x^k)\\
	&+ 	\beta_k\left[f(\u) + (n\beta_k)^{\gamma-1} D_H(\u,\z^k)
	- (n\beta_k)^{\gamma-1} D_H(\u,\tilde{\z}^{k+1})
	\right],
	\end{align*}
	where $(i)$ is because the $i_k$-th coordinate is selected uniformly at random, $(ii)$ is due to the convexity of $f$, $(iii)$ is due to the definition of $\tilde{z}^{k+1}$ in E.q. \eqref{eq:zhat^+} and applying Lemma~\ref{lemma:prox} with $\phi(\u) = (n\beta_k)^{1-\gamma}\left[f(\y^k) + \inn{\nabla f(\y^k)}{\u - \y^k} + \delta_{\mathcal{X}}(\u)\right]$, and $(iv)$ is due to the convexity of $f$. Subtracting $f(\u)$ on both sides gives us
	\begin{align*}
	\E_{i_k} f(\x^{k+1}) -& f(\u) \leq (1-\beta_k)( f(\x^k) - f(\u))\\
	&+ n^{\gamma-1}\beta_k^{\gamma} D_H(\u,\z^k)
	- n^{\gamma-1}\beta_k^{\gamma} D_H(\u,\tilde{\z}^{k+1}).
	\end{align*}
	Dividing $\beta_k^\gamma$ on both sides, we have
	\begin{align}
	&\frac{1}{\beta_k^{\gamma}}\E_{i_k} \left[f(\x^{k+1}) - f(\u)\right] \leq 
	\frac{1-\beta_k}{\beta_k^{\gamma}}( f(\x^k) - f(\u))\nonumber\\
	&+ n^{\gamma-1} D_H(\u,\z^k)
	- n^{\gamma-1} D_H(\u,\tilde{\z}^{k+1}).
	\label{eq:f(x^{k+1})}
	\end{align}
	
	Taking the expectation of $D_H(\u, \z^{k+1})$ with respect to $\E_{i_k}$ yields
	\begin{align*}
	&\E_{i_k} [D_H(\u, \z^{k+1})]  \\
	=& \E_{i_k} \left[ D_H(\u,\z^k) - L_{i_k} D_h(\u_{i_k}, \z^k_{i_k}) + L_{i_k} D_h(\u_{i_k}, \tilde{\z}^{k+1}_{i_k})  \right] \\
	=&\sum_{i_k=1}^n \frac{1}{n} \left[ D_H(\u,\z^k) - L_{i_k} D_h(\u_{i_k}, \z^k_{i_k}) + L_{i_k} D_h(\u_{i_k}, \tilde{\z}^{k+1}_{i_k} )  \right]\\
	=& \frac{1}{n} \left[ n D_H(\u,\z^k) - D_H(\u, \z^k) + D_H(\u, \tilde{\z}^{k+1})  \right]\\
	=& D_H(\u,\z^k) - \frac{1}{n} \left[ D_H(\u, \z^k) - D_H(\u, \tilde{\z}^{k+1})  \right] .
	\end{align*}
	Multiplying both sides by $n^\gamma$, we obtain
	\begin{align}
	&n^\gamma\E_{i_k} [D_H(\u, \z^{k+1})]\nonumber\\
	=& n^\gamma D_H(\u,\z^k) - n^{\gamma-1}\left[ D_H(\u, \z^k) - D_H(\u, \tilde{\z}^{k+1})  \right]
	\label{eq:D_H}
	\end{align}
	
	Combining \eqref{eq:D_H} with \eqref{eq:f(x^{k+1})}, we have
	\begin{align}
	\E_{i_k} \left[\frac{1}{\beta_k^{\gamma}}(f(\x^{k+1}) - f(\u)) + n^\gamma D_H(\u, \z^{k+1}) \right]\nonumber\\
	\leq \frac{1-\beta_k}{\beta_k^{\gamma}}( f(\x^k) - f(\u)) + n^\gamma D_H(\u,\z^k)
	\end{align}
	Finally applying the condition in Step 4 of Algorithm~\ref{alg:arbcd} yields the desired result.
\end{proof}

\subsection{Proof of Theorem~\ref{thm:acc}}
\begin{proof}
	Taking the expectation with respect to $\{i_0, i_1,\cdots,\}$ yields
	\begin{align}
	\E \left[\frac{1-\beta_{k+1}}{\beta_{k+1}^{\gamma}}(f(\x^{k+1}) - f(\u)) + n^\gamma D_H(\u, \z^{k+1}) \right] \nonumber\\
	\leq \E\left[\frac{1-\beta_k}{\beta_k^{\gamma}}( f(\x^k) - f(\u)) + n^\gamma D_H(\u,\z^k) \right].
	\label{eq:rec}
	\end{align}
	The direct consequence of E.q. \eqref{eq:rec} is, for any $\u$,
	\begin{align*}
	\E \left[\frac{1-\beta_{k+1}}{\beta_{k+1}^{\gamma}}(f(\x^{k+1}) - f(\u)) + n^\gamma D_H(\u, \z^{k}) \right] \nonumber \\
	\leq \frac{1-\beta_0}{\beta_0^{\gamma}}( f(\x^0) - f(\u)) + n^\gamma D_H(\u,\z^0) .
	\end{align*}
	Using $D_H(\u, \z^{k+1})\geq 0$, and the initialization $\beta_0 = 1$ and $\z^0 = \x^0$, we obtain
	\begin{align*}
	\E \left[\frac{1-\beta_{k+1}}{\beta_{k+1}^{\gamma}}(f(\x^{k+1}) - f(\u)) \right] \leq n^\gamma D_H(\u,\x^0),
	\end{align*}
	which implies
	\begin{align*}
	\E \left[f(\x^{k+1}) - f(\u) \right] \leq n^\gamma \beta_{k}^\gamma D_H(\u,\x^0) \\
	= \left(\frac{n\gamma}{k+\gamma}\right)^\gamma D_H(\u,\x^0).
	\end{align*}
\end{proof}

\subsection{Proof of Proposition~\ref{prop:eff}}
\begin{proof}	
	It is straightforward to see that $\x^0 = \y^0 = \z^0 = \v^0$. Suppose the recursive hypotheses hold for the $k$-th iteration. From the optimality of E.q. \eqref{eq:d^k}, we have
	\begin{align}
	&\inn{\nabla_{i_k} f(\beta_k^\gamma \u^k + \v^k)}{\d^k_{i_k}} + (n\beta_k)^{\gamma-1}L_{i_k}D_h(\v^k_{i_k}+\d^k_{i_k},\v_{i_k}^k) \nonumber\\
	\overset{(i)}{\leq} &\inn{\nabla_{i_k} f(\beta_k^\gamma \u^k + \v^k)}{\z^{k+1}_{i_k} - \z^k_{i_k}} + (n\beta_k)^{\gamma-1}L_{i_k}D_h(\z^{k+1}_{i_k},\v_{i_k}^k)\nonumber\\
	\overset{(ii)}{=}	&\inn{\nabla_{i_k} f(\y^k)}{\z^{k+1}_{i_k} - \z^k_{i_k}} + (n\beta_k)^{\gamma-1}L_{i_k}D_h(\z^{k+1}_{i_k},\z_{i_k}^k),
	\label{eq:left}
	\end{align}
	where $(i)$ is due to the optimality, and $(ii)$ is due to the recursive hypotheses. Similarly, from the optimality of E.q. \eqref{eq:z^+}, we obtain
	\begin{align}
	&\inn{\nabla_{i_k} f(\y^k)}{\z^{k+1}_{i_k} - \z^k_{i_k}} + (n\beta_k)^{\gamma-1}L_{i_k}D_h(\z^{k+1}_{i_k},\z_{i_k}^k)\nonumber\\
	\overset{(i)}{\leq}& \inn{\nabla_{i_k} f(\y^k)}{\d^k_{i_k}} + (n\beta_k)^{\gamma-1}L_{i_k}D_h(\z^{k}_{i_k} + d^k_{i_k},\z_{i_k}^k)\nonumber\\
	\overset{(ii)}{=}& \inn{\nabla_{i_k} f(\beta_k^\gamma \u^k + \v^k)}{\d^k_{i_k}} + (n\beta_k)^{\gamma-1}L_{i_k}D_h(\v^k_{i_k}+\d^k_{i_k},\v_{i_k}^k), \label{eq:right}
	\end{align}
	where $(i)$ is due to the optimality, and $(ii)$ is due to the recursive hypotheses. Combing \eqref{eq:left} and \eqref{eq:right} yields
	\begin{align*}
	\z_{i_k}^{k+1} = \z_{i_k}^k + \d_{i_k}^k = \v_{i_k}^k + \d_{i_k}^k = \v_{i_k}^{k+1},
	\end{align*}
	or equivalently
	\begin{align*}
	\z^{k+1} = \v^{k+1}.
	\end{align*}
	From Step 3 of Algorithm~\ref{alg:efficient}, we have
	\begin{align}
	\u^{k+1} = \u^k - \frac{1-n\beta_{k}}{\beta_{k}^\gamma}(\v^{k+1} - \v^k). \label{eq:u^+}
	\end{align}
	Then, we have
	\begin{align*}
	\beta_{k}^\gamma \u^{k+1} + \v^{k+1}
	\overset{(i)}{=}&\beta_{k}^{\gamma}\left(\u^k - \frac{1-n\beta_{k}}{\beta_{k}^\gamma}(\v^{k+1} - \v^k) \right) + \v^{k+1}\\
	=&\beta_{k}^{\gamma} \u^k - (1-n\beta_{k} ) (\v^{k+1} - \v^k) + \v^{k+1}\\
	=&\beta_{k}^{\gamma} \u^k + \v^k + n\beta_{k}  (\v^{k+1} -  \v^{k}) \\
	\overset{(ii)}{=}& \y^k + n\beta_{k}  (\z^{k+1} -  \z^{k}) \\
	=&\x^{k+1},
	\end{align*}
	where $(i)$ is due to E.q. \eqref{eq:u^+} and $(ii)$ is due to the recursive hypotheses.
	
	Finally, we have
	\begin{align*}
	\beta_{k+1}^\gamma\u^{k+1} + \v^{k+1}\overset{(i)}{=}&\frac{\beta_{k+1}^\gamma}{\beta_{k}^\gamma}(\x^{k+1} - \v^{k+1}) + \v^{k+1}\\
	\overset{(ii)}{=}&(1-\beta_{k+1})(\x^{k+1} - \v^{k+1}) + \v^{k+1}\\
	=&(1-\beta_{k+1}) \x^{k+1}  + \beta_{k+1}\v^{k+1}\\
	\overset{(iii)}{=}&(1-\beta_{k+1}) \x^{k+1}  + \beta_{k+1}\z^{k+1}\\
	=&\y^{k+1},
	\end{align*}
	where $(i)$ and $(iii)$ is due to recursive hypotheses, and $(ii)$ is due to Step 4 of Algorithm~\ref{alg:efficient}.
\end{proof}

\subsection{Proof of Lemma~\ref{lemma:relative smooth}}
\begin{proof}
	Let us fix coordinate $j$. Then we have
	\begin{align*}
		f_j(\x_j) = \sum_{i=1}^M \left(\b_i \log\left(\frac{\b_i}{\inn{\a_i}{\x}}\right) + \inn{\a_i}{\x} - \b_i\right),
	\end{align*}
	where $\a_i$ is the $i$-th row of $\A$. The first- and second-order derivatives of $f_j$ are given by
	\begin{align*}
		f_j^\prime(\x_j) &=\sum_{i=1}^M \left(1-\frac{\b_i}{\inn{\a_i}{\x}}\right)\a_{ij},\\
		f_j^{\prime\prime}(\x_j) &=\sum_{i=1}^M \frac{\b_i \a_{ij}^2}{\inn{\a_i}{\x}^2}.
	\end{align*}
	It follows from the nonnegativity of $\A$ and $\x$ that we have
	\begin{align*}
		\frac{ \a_{ij}^2}{\inn{\a_i}{\x}^2}\leq \frac{1}{\x^2_{j} }.
	\end{align*}
	Applying the inequality above yields
	\begin{align*}
		f_j^{\prime\prime}(\x_j) =\sum_{i=1}^M \frac{\b_i \a_{ij}^2}{\inn{\a_i}{\x}^2}
		\leq \left(\sum_{i=1}^M  \b_i\right) \frac{1}{\x^2_{j} }
		= \left(\sum_{i=1}^M  \b_i\right) h_j^{\prime\prime} (\x_j).
	\end{align*}
\end{proof}

\subsection{Proof of Lemma~\ref{lemma:relative entropy}}
\begin{proof}
	Fixing the $j$-th coordinate of $\x$, define $f_j(\x_j)$ as follows
	\begin{align*}
	f_j(\x_j) = \sum_{i=1}^M \left(\inn{\a_i}{\x} \log\left(\frac{\inn{\a_i}{\x}}{\b_i}\right) +\b_i - \inn{\a_i}{\x} \right.
	\end{align*}
	Then the first- and second-derivatives of $f_j$ are given by
	\begin{align*}
	f_j^\prime(\x_j) &=\sum_{i=1}^M \a_{ij}\left(\log\frac{\inn{\a_i}{\x}}{\b_i}\right),\\
	f_j^{\prime\prime}(\x_j) &=\sum_{i=1}^M \frac{\a_{ij}^2}{\inn{\a_i}{\x}}.
	\end{align*}
	Using the nonnegativity of $\A$ and $\x$, we obtain $\a_{ij} \x_j\leq \inn{\a_i}{\x}$, which further implies
	\begin{align}
		\frac{\a^2_{ij}}{\inn{\a_i}{\x}}\leq\frac{\a_{ij}}{\x_j}.
	\end{align}
	Invoking the inequality above, we obtain the desired result
	\begin{align*}
	f_j^{\prime\prime}(\x_j) =\sum_{i=1}^M \frac{ \a_{ij}^2}{\inn{\a_i}{\x}}
	\leq \sum_{i=1}^M \frac{\a_{ij}}{\x_j}
	=\left(\sum_{i=1}^M\a_{ij}\right) h_j^{\prime\prime} (\x_j).
	\end{align*}
\end{proof}

\bibliographystyle{IEEE-unsorted}
\bibliography{refs}
\end{document}